\documentclass[a4paper,10.5pt]{article}
\usepackage[pagewise]{lineno}
\usepackage{amssymb}
\usepackage{latexsym, bm}
\usepackage{multicol}
\usepackage{indentfirst}
\usepackage{amssymb,amsfonts}
\usepackage{amsmath}
\usepackage{setspace}

\usepackage{amssymb}
\title{A regularity criterion for three-dimensional micropolar fluid equations in Besov spaces of negative regular indices}
\author{Maria Alessandra Ragusa\footnote{maragusa@dmi.unict.it}\\
{\small Department of Mathematics, University of Catania,}\\
{\small Viale Andrea Doria No. 6, Catania 95128, Italy}\\
{\small RUDN University,}\\
{\small6 Miklukho -Maklay St, Moscow, 117198, Russia}
\and
Fan Wu\footnote{wufan0319@yeah.net}\\
{\small School of Mathematics and Statistics, Hunan Normal University, }\\
{\small Changsha, Hunan 410081, China }}

\date{}

\begin{document}
\maketitle
{\bf Abstract:}
In this article, we study regularity criteria for the 3D  micropolar fluid equations in terms of one partial derivative of the velocity. It is proved that if
\begin{equation*}
\int^{T}_{0}\|\partial_{3}u\|^{\frac{2}{1-r}}_{\dot{B}^{-r}_{\infty,\infty}} dt<\infty \quad \text{with} \quad  0< r<1,
\end{equation*}
then, the solutions of the micropolar fluid equations actually are smooth on $(0, T)$. This improves and extends many previous results.
\medskip

\medskip

\section{Introduction}
We are interested in the regularity of weak solutions to the incompressible micropolar fluid equations in $\mathbb{R}^3$:
\begin{spacing}{1.5}
\begin{equation}
\left\{
             \begin{array}{lr}
             \partial_t u + (u\cdot\nabla) u +\nabla \pi
 =  \Delta u +  \nabla \times \omega,&\\
\partial_t \omega + (u\cdot \nabla) \omega + 2 \omega
 =  \Delta \omega+  \nabla \nabla\cdot \omega +  \nabla \times u,&\\
             \nabla\cdot u=0,&\\
             u(x,0)=u_{0}(x),\omega(x,0)=\omega_{0}(x),&
\end{array}
\right.\tag{1.1}
\end{equation}
\end{spacing}
where, for $\mathbf{x}\in \mathbb{R}^3$ and $t\ge 0$,
$u=u(\mathbf{x},t), \omega=\omega(\mathbf{x},t)$  and $\pi=\pi(\mathbf{x},t)$
denote the velocity field, the micro-rotation field and the pressure,
respectively.
\par Micropolar fluids represent a class of fluids with nonsymmetric stress tensor (called polar fluids) such as fluids consisting of suspending particles, dumbbell molecules, (see, e.g.,
\cite{ Erdo, Er, Er2}). Mathematically, many authors \cite{RM,YK,DZ} treated the well-posedness
and large-time behaviour of solutions to system (1.1). However, the issue of global regularity
of weak solutions to (1.1) remains an open problem. In the absence of micro-rotational
effects ($\omega=0$), this system reduces to well-known incompressible Navier Stokes equations. Therefore, it is an interesting thing that regularity of a given weak solution of the 3D micropolar fluids or the 3D Navier-Stokes equations can be shown under some additional conditions, and over the years different criteria for regularity of the weak solutions have been proposed. Some fundamental Serrin-type regularity criteria for system $(1.1)$ in terms of the velocity only were carried out in \cite{YB,YB1,DZ1} independently. Recently, some
improvement and extension were made on the basis of the present paper (see e.g. \cite{DZ2,GS,ZZ,WY} and the references therein) were derived to guarantee the regularity of the weak solution.
\par In this paper, we are concerned with the regularity conditions of partial components of three-dimensional micropolar fluid equations. In this respect, Jia et al. \cite{JZZ}
showed that if
\begin{equation}
\partial_3u\in L^{\frac{2p}{p-3}}(0,T;L^{p,\infty}) \quad  \text{with} \quad 3<p\leq\infty ,\tag{1.2}
\end{equation}
or
\begin{equation}
\partial_1u_1, \partial_2u_2\in L^{\frac{2p}{2p-3}}(0,T;L^{p,\infty}) \quad \text{with}\quad \frac{3}{2}<p\leq\infty,\tag{1.3}
\end{equation}
then the solution is smooth on (0, T). Later, in \cite{GS1}, the authors extended the regularity criterion $(1.2)$ to the Morrey-Campanato spaces. More precisely, they obtained following interesting regularity criterion
 \begin{equation}
\partial_3u\in L^{\frac{2}{1-r}}(0,T;\dot{M}_{2,\frac{3}{r}}) \quad \text{with} \quad  0< r<1.\tag{1.4}
\end{equation}
Recently, many authors refine earlier results from the Lebesgue
spaces to Besov spaces. Zhang \cite{ZZ} and Yuan \cite{YB2} refined the BKM criterion
\begin{equation}
u\in L^{\frac{2}{1+r}}(0,T;\dot{B}^{r}_{\infty,\infty}) \quad \text{with} \quad 0<r<1,\tag{1.5}
\end{equation}
and
 \begin{equation}
\nabla_{h}u,\nabla_{h}\omega\in L^{1}(0,T;\dot{B}^{0}_{\infty,\infty}), \tag{1.6}
\end{equation}
or
\begin{equation}
\nabla_{h}u,\nabla_{h}\omega\in L^{\frac{8}{3}}(0,T;\dot{B}^{-1}_{\infty,\infty}), \tag{1.7}
\end{equation}
respectively. Readers interested in regularity criteria for other incompressible models may refer to \cite{GLR,GLR1,GLR2,LQ,GS,CW,GS2,JL,JX,ZY}.
\par Motivated by the references mentioned above, the purpose of the present paper
is to extend the regularity criterion of weak solutions to the micropolar fluid equations $(1.1)$ in terms of one partial derivative of the velocity in the framework of the homogeneous Besov space, which improves the recently results of \cite{JZZ,GS1}.
\medskip
\par Now, we state our results as follows.
\medskip
\\\textbf{Theorem 1.1}\quad Assume that the initial velocity and micro-rotation field $(u_{0},\omega_0)\in H^{1}(R^{3})$ with $\nabla\cdot u_{0}=0$. Let $(u,\omega)$ be a weak solution to the system $(1.1)$ on some interval $[0,T]$ with $0<T< \infty$. If additionally
\begin{equation}
\int^{T}_{0}\|\partial_{3}u\|^{\frac{2}{1-r}}_{\dot{B}^{-r}_{\infty,\infty}} dt<\infty \quad \text{with} \quad  0< r<1\tag{1.8}
\end{equation}
then, the solution $(u,\omega)$ smooth on $[0,T]$.
\medskip
\\\textbf{Remark 1.1}\quad For $2\leq p\leq\frac{3}{r}$ and $0\leq r<\frac{3}{2}$, we have following inclusion relations:
\[L^{\frac{3}{r}}(\mathbb{R}^3)\hookrightarrow L^{\frac{3}{r},\infty}(\mathbb{R}^3)\hookrightarrow\dot{M}_{p,\frac{3}{r}}(\mathbb{R}^3)\hookrightarrow\dot{X}_r(\mathbb{R}^3)\hookrightarrow\dot{M}_{2,\frac{3}{r}}(\mathbb{R}^3)\hookrightarrow\dot{B}^{-r}_{\infty,\infty}(\mathbb{R}^3).\]
Therefore, Theorem 1.1 can be regarded as an extension of condition $(1.2)$ and $(1.4)$.
 \medskip
\\\textbf{Remark 1.2}\quad  If let $\omega= 0$ in system $(1.1)$, it is clear that Theorem 1.1 improve and extend the results in \cite{GS,LQ} for the 3D incompressible Navier-Stokes equations.
 \medskip
\\\textbf{Remark 1.3}\quad  Different from the literatures \cite{JZZ,GS1}, we estimate the first equation of $(1.1)$ to obtain the priori estimate of the velocity field $u$, and then use the iterative method to estimate the micro-rotation field $\omega$ to obtain the regularity condition of the solution for the three-dimensional micropolar fluid equation.
\medskip
\par Next, in order to derive the regularity criterion of weak solutions to the micropolar
fluid equations $(1.1)$, we introduce the definition of a weak solution.
\medskip
\\\textbf{Definition 1.1}\quad Let $u_0\in L^2(\mathbb{R}^3)$ with $\nabla\cdot u_0=0$, $\omega_0=L^2(\mathbb{R}^3)$. A measurable $\mathbb{R}^3$-valued pair $(u,\omega)$ is called a weak solution to system $(1.1)$ on $(0,T)$, provided that following two conditions hold,
\medskip
\par $(1)$\quad $(u,\omega)\in L^{\infty}(0,T;L^2(\mathbb{R}^3))\cap L^{2}(0,T;H^1(\mathbb{R}^3))$;
\medskip
\par $(2)$\quad $(u,\omega)$ verifies system $(1.1)$ in the distributions sense.
\medskip
\par The rest of this paper is organized as follows. We shall present some preliminary
on functional settings and some useful lemmas in Section 2, while the proof of Theorem 1.1 will be presented in Section 3.
\section{Preliminaries}
In this section, we give some lemmas and introduce some basic facts on Littlewood-Paley theory, which will be used in the proof of our main results.
\par \par Let $\mathcal{S}(\mathbb{R}^3)$ be the Schwartz class of rapidly decreasing functions. Given $f\in \mathcal{S}(\mathbb{R}^3)$, its Fourier transform $\mathcal{F}f=\hat{f}$ is defined as
\[\hat{f}(\xi)=\int_{\mathbb{R}^3}f(x)e^{-ix\cdot\xi}dx.\]
Let $(\chi,\varphi)$ be a couple of smooth functions valued in $[0,1]$ such that $\chi$ is supported in $B=\{\xi \in \mathbb{R}^3:|\xi|\leq\frac{4}{3}\}$,
$\varphi$ is supported in $\mathcal{C}=\{\xi \in \mathbb{R}^3:\frac{3}{4}\leq|\xi|\leq\frac{8}{3}\}$ such that
\[\chi(\xi)+\sum\limits_{j\geq 0}\varphi(2^{-j}\xi)=1,\quad \forall\xi\in \mathbb{R}^3,\]
\[\sum\limits_{j\in Z}\varphi(2^{-j}\xi)=1, \quad \forall\xi\in \mathbb{R}^3\setminus \{0\}.\]
Denoting $\varphi_{j}=\varphi(2^{-j}\xi)$, $h=\mathcal{F}^{-1}\varphi$ and $\widetilde{h}=\mathcal{F}^{-1}\chi$, the
dyadic blocks are defined as follows, respectively.
\[\dot{\Delta}_{j}f=\varphi(2^{-j}D)f=2^{3j}\int_{\mathbb{R}^3}h(2^{j}y)f(x-y)dy,\]
\[\dot{S}_{j}f=\sum \limits_{k\leq j-1}\Delta_{j}f=\chi(2^{-j}D)f=2^{3j}\int_{\mathbb{R}^3}\widetilde{h}(2^{j}y)f(x-y)dy.\]
Furthermore, the above dyadic decomposition has nice properties of quasi-orthogonality,
namely,
\begin{center}
$\dot{\Delta}_{j}\dot{\Delta}_{q}f\equiv0$ \quad with $|j-q|\geq2.$
\end{center}
We have the following formal decomposition:
\[f=\sum \limits^{\infty}_{-\infty}\dot{\Delta}_{j}f\quad for \quad f\in\mathcal{S'}(\mathbb{R}^3)\setminus \mathcal{P}(\mathbb{R}^3),\]
where $\mathcal{P}(\mathbb{R}^3)$ is the set of polynomials, which can be found in \cite{BC} about the details of the Littlewood-Paley decomposition theory.
\par Lets now recall the homogeneous Besov spaces.
\\\\\textbf{Definition 2.1}\quad Let $s\in \mathbb{R}, (p,q)\in [1,\infty]^{2}$, the homogeneous Besov space $\dot{B}_{p,q}^{s}$ is defined by
\[\dot{B}_{p,q}^{s}=\{f\in \mathcal{Z^{'}}(\mathbb{R}^3);\|f\|_{\dot{B}_{p,q}^{s}}<\infty\},\]
where
$$\|f\|_{\dot{B}_{p,q}^{s}}=
\begin{cases}
(\sum\limits_{j\in Z}2^{jsq}\|\Delta_{j}f\|_{p}^{q})^{\frac{1}{q}}& \text{$q<\infty$}, \\
\sup\limits_{j\in Z}2^{js}\|\Delta_{j}f\|_{p}& \text{$q=\infty$.}
\end{cases}$$
Here $ \mathcal{Z^{'}}(\mathbb{R}^3)$ is the dual space of
\[\mathcal{Z}(\mathbb{R}^3)=\{f\in\mathcal{S}(\mathbb{R}^3);D^{\alpha}\widehat{f}(0)=0, \forall\alpha\in \mathbb{N}^3\} .\]
Finally, for completeness, we give following Bernstein's inequalities and interpolation inequality (see \cite{BC}), which plays an important role in the proof of main results.
\medskip
\\\textbf{Lemma 2.1}\quad Let $1\leq q<p<\infty$ and $\alpha$ be a positive real number. A constant $C$ exists such that
\begin{equation}
\|f\|_{L^{p}}\leq C\|f\|^{1-\theta}_{\dot{B}^{-\alpha}_{\infty,\infty}}\|f\|^{\theta}_{\dot{B}^{\beta}_{q,q}},\tag{2.1}
\end{equation}
where $\beta=\alpha(\frac{p}{q}-1)$ and $\theta=\frac{q}{p}$.
\medskip
\\\textbf{Lemma 2.2} \cite{CW}\quad Let $\mu,\theta,\lambda$ and $\kappa$ be four numbers satisfying the following relations
\begin{displaymath}
1\leq\mu,\theta,\lambda,\kappa<\infty,\quad\frac{1}{\theta}+\frac{1}{\lambda}+\frac{1}{\kappa}>1
\quad \text{and} \quad 1+\frac{3}{\mu}=\frac{1}{\theta}+\frac{1}{\lambda}+\frac{1}{\kappa}.
\end{displaymath}
Assume that $\varphi(x)=\varphi(x_{1},x_{2},x_{3})$ with $\partial_{1}\varphi\in L^{\theta}(\mathbb{R}^3),\partial_{2}\varphi\in L^{\lambda}(\mathbb{R}^3)$ and $\partial_{3}\varphi\in L^{\kappa}(\mathbb{R}^3)$.
Then, there exists a constant $C=C(\theta,\lambda,\kappa)$ such that
\begin{equation}
\|\varphi\|_{L^{\mu}}\leq C\|\partial_{1}\varphi\|^{\frac{1}{3}}_{L^{\theta}}
\|\partial_{2}\varphi\|^{\frac{1}{3}}_{L^{\lambda}}
\|\partial_{3}\varphi\|^{\frac{1}{3}}_{L^{\kappa}}.\tag{2.2}
\end{equation}
Especially, when $\theta=\lambda=2$ and $1\leq\kappa<\infty$, there exists a constant $C=C(\kappa)$ such that
\begin{equation}
\|\varphi\|_{L^{3\kappa}}\leq C\|\partial_{1}\varphi\|^{\frac{1}{3}}_{L^{2}}
\|\partial_{2}\varphi\|^{\frac{1}{3}}_{L^{2}}
\|\partial_{3}\varphi\|^{\frac{1}{3}}_{L^{\kappa}},\tag{2.3}
\end{equation}
which holds for any $\varphi$ with $\partial_{1}\varphi\in L^{2}(\mathbb{R}^3), \partial_{2}\varphi\in L^{2}(\mathbb{R}^3)$ and $\partial_{3}\varphi\in L^{\kappa}(\mathbb{R}^3)$.
\section{Proof of Theorem 1.1}
This section is devoted to the proof of Theorem 1.1. It is based on the establishment of a priori estimates under condition $(1.8)$.
\par Multiplying the first equation of $(1.1)$ by $u$, and the second equation of $(1.1)$ by $\omega$, integrating over $\mathbb{R}^3$, then we add the resulting equations, it yields
\begin{equation}
\frac{1}{2}\frac{d}{dt}(\|u\|^{2}_{L^{2}}+\|\omega\|^{2}_{L^{2}})+\|\nabla u\|^{2}_{L^{2}}+\|\nabla \omega\|^{2}_{L^{2}}+2\|\omega\|^{2}_{L^{2}}+\|\nabla\cdot\omega\|^{2}_{L^{2}}=0.\tag{3.1}
\end{equation}
Integrating in time, we get the fundamental energy estimate
 \begin{equation}
\frac{1}{2}\|(u,\omega)(t)\|^{2}_{L^{2}}+\int^{t}_{0}\|(\nabla u,\nabla \omega)(s)\|^{2}_{L^{2}}ds+\int^{t}_{0}2\|\omega(s)\|^{2}_{L^{2}}ds
+\int^{t}_{0}\|\nabla\cdot\omega(s)\|^{2}_{L^{2}}ds\leq\frac{1}{2}\|(u_{0},\omega_{0})\|^{2}_{L^{2}}.\tag{3.2}
\end{equation}
\medskip
\\\textbf{Lemma 3.1}\quad Assume that $u_0,\omega_0\in H^1(\mathbb{R}^3)$ with $\nabla\cdot u_0$, let $T>0$, $(u,\omega)$ is the weak solution to the system $(1.1)$. If there exists a constant $K>0$ such that the following conditions holds,
\begin{equation}
\int^{T}_{0}\|\partial_{3}u\|^{\frac{2}{1-r}}_{\dot{B}^{-r}_{\infty,\infty}} dt<K \quad \text{with} \quad  0< r<1,\tag{3.3}
\end{equation}
then, we have
 \begin{equation}
\begin{split}
\sup\limits_{0\leq t\leq T}\|\nabla u\|^{2}_{L^{2}}+\int^T_0\|\Delta u\|^2_{L^2}dt\leq C,
\end{split}\tag{3.4}
\end{equation}
where constant $C$ depending only on $u_0$, $\omega_0$, $T$ and $K$.\\
\\\textbf{Proof:}
Differentiating the first equation of system $(1.1)$ about space variable $x_{3}$, then multiplying the resulting equation by $\partial_{3}u$, and integrating it we get
\begin{equation}
\begin{split}
\frac{1}{2}\frac{d}{dt}\|\partial_{3} u\|^{2}_{L^{2}}+\|\nabla \partial_{3}u\|^{2}_{L^{2}}&=-\int_{\mathbb{R}^3}\partial_{3}(u\cdot\nabla u)\cdot\partial_{3} u dx+\int_{\mathbb{R}^3}\partial_{3}(\nabla\times\omega)\cdot\partial_{3} u dx\\
&=I_{1}+I_{2}.
\end{split}\tag{3.5}
\end{equation}
By Young inequality, standard multiplicative inequality, Lemma 2.1 and Lemma 2.2, we can estimate $I_{1}$ and $I_{2}$ as follows:
\begin{equation}
\begin{split}
|I_{1}|&=\int_{\mathbb{R}^3}\partial_{3}u\cdot\nabla u\cdot\partial_{3} u dx\\
&\leq\|\partial_{3}u\|^{2}_{L^{4}}\|\nabla u\|_{L^{2}}\\
&\leq\|\partial_{3}u\|_{\dot{B}^{-r}_{\infty,\infty}}\|\partial_{3}u\|_{\dot{H}^{r}}\|\nabla u\|_{L^{2}}\\
&\leq\|\partial_{3}u\|_{\dot{B}^{-r}_{\infty,\infty}}\|\partial_{3}u\|^{1-r}_{L^{2}}
\|\nabla\partial_{3}u\|^{r}_{L^{2}}\|\nabla u\|_{L^{2}}\\
&\leq C(\|\nabla\partial_{3}u\|^{2}_{L^{2}})^{\frac{r}{2}}(\|\partial_{3}u\|^{\frac{2}{2-r}}_{\dot{B}^{-r}_{\infty,\infty}}
\|\partial_{3}u\|^{2\cdot\frac{1-r}{2-r}}_{L^{2}}\|\nabla u\|^{\frac{2}{2-r}}_{L^{2}})^{\frac{2-r}{2}}\\
&\leq \frac{1}{4}\|\nabla\partial_{3}u\|^{2}_{L^{2}}+C\|\partial_{3}u\|^{\frac{2}{2-r}}_{\dot{B}^{-r}_{\infty,\infty}}
\|\partial_{3}u\|^{2\cdot\frac{1-r}{2-r}}_{L^{2}}\|\nabla u\|^{\frac{2}{2-r}}_{L^{2}}\\
&\leq \frac{1}{4}\|\nabla\partial_{3}u\|^{2}_{L^{2}}+C[(\|\partial_{3}u\|^{\frac{2}{1-r}}_{\dot{B}^{-r}_{\infty,\infty}})^{\frac{1-r}{2-r}}
(\|\nabla u\|^{2}_{L^{2}})^{\frac{1}{2-r}}]
\|\partial_{3}u\|^{2\cdot\frac{1-r}{2-r}}_{L^{2}}\\
&\leq \frac{1}{4}\|\nabla\partial_{3}u\|^{2}_{L^{2}}+C(1+\|\partial_{3}u\|^{2}_{L^{2}})(\|\partial_{3}u\|^{\frac{2}{1-r}}_{\dot{B}^{-r}_{\infty,\infty}}
+\|\nabla u\|^{2}_{L^{2}}),
\end{split}\tag{3.6}
\end{equation}
and
\begin{equation}
\begin{split}
|I_{2}|=|\int_{\mathbb{R}^3}\partial_{3}(\nabla\times\omega)\cdot\partial_{3} u dx|\leq \|\partial_{33}u\|_{L^{2}}\|\nabla\omega\|_{L^{2}}\leq \frac{1}{4}\|\nabla\partial_{3} u\|^{2}_{L^{2}}+C\|\nabla\omega\|^{2}_{L^{2}}.
\end{split}\tag{3.7}
\end{equation}
Combining $(3.5)$, $(3.6)$ and $(3.7)$ together, then we deduce that
\begin{equation}
\begin{split}
\frac{d}{dt}(1+\|\partial_{3} u\|^{2}_{L^{2}})+\|\nabla \partial_{3}u\|^{2}_{L^{2}}&\leq C(1+\|\partial_{3} u\|^{2}_{L^{2}})(1+\|\nabla u\|^{2}_{L^{2}}+\|\nabla\omega\|^{2}_{L^{2}}+\|\partial_{3}u\|^{\frac{2}{1-r}}_{\dot{B}^{-r}_{\infty,\infty}}).
\end{split}\tag{3.8}
\end{equation}
Applying the Gronwall inequality to $(3.8)$ yields that
\begin{equation}
\begin{split}
(1+\|\partial_{3} u\|^{2}_{L^{2}})&\leq (1+\|\partial_{3} u_{0}\|^{2}_{L^{2}})\exp\{C\int^{T}_{0}(1+\|\nabla u\|^{2}_{L^{2}}+\|\nabla\omega\|^{2}_{L^{2}}+\|\partial_{3}u\|^{\frac{2}{1-r}}_{\dot{B}^{-r}_{\infty,\infty}})d\tau\},
\end{split}\tag{3.9}
\end{equation}
which implies that
\begin{equation}
\begin{split}
\sup\limits_{0\leq t\leq T}\|\partial_{3} u(t)\|^{2}_{L^{2}}+\int^{T}_{0}\|\nabla\partial_{3} u(t)\|^{2}_{L^{2}}dt&\leq C.
\end{split}\tag{3.10}
\end{equation}
Next, multiplying the first equation of system $(1.1)$ by $\Delta u$, integrating over $\mathbb{R}^3$, we get
 \begin{equation}
 \begin{split}
\frac{1}{2}\frac{d}{dt}\|\nabla u\|^{2}_{L^{2}}+\|\Delta u\|^{2}_{L^{2}}&=\int_{\mathbb{R}^3}u\cdot\nabla u\cdot\Delta u dx-\int_{\mathbb{R}^3}(\nabla\times\omega)\cdot\Delta u dx\\
&=J_{1}+J_{2}.
\end{split}\tag{3.11}
\end{equation}
For $J_{1}$, by using the interpolation inequality and Lemma 2.2 (with $\kappa=2$ in (2.3)), it is not difficult to see that
\begin{equation}
\begin{split}
J_{1}&\leq C\|\nabla u\|^{3}_{L^{3}}\leq C\|\nabla u\|^{\frac{3}{2}}_{L^{2}}\|\nabla u\|^{\frac{3}{2}}_{L^{6}}\\
&\leq C\|\nabla u\|^{\frac{3}{2}}_{L^{2}}\|\nabla\partial_{1}u\|^{\frac{1}{2}}_{L^{2}}
\|\nabla\partial_{2}u\|^{\frac{1}{2}}_{L^{2}}\|\nabla\partial_{3}u\|^{\frac{1}{2}}_{L^{2}}\\
&\leq C\|\nabla u\|^{\frac{3}{2}}_{L^{2}}\|\nabla^{2}u\|_{L^{2}}
\|\nabla\partial_{3}u\|^{\frac{1}{2}}_{L^{2}}\\
&=C(\|\nabla^{2}u\|^{2}_{L^{2}})^{\frac{1}{2}}(\|\nabla u\|^{3}_{L^{2}}\|\nabla\partial_{3}u\|_{L^{2}})^{\frac{1}{2}}\\
&\leq\frac{1}{4}\|\Delta u\|^{2}_{L^{2}}+C\|\nabla u\|^{3}_{L^{2}}\|\nabla\partial_{3}u\|_{L^{2}}\\
&\leq\frac{1}{4}\|\Delta u\|^{2}_{L^{2}}+C\|\nabla u\|^{2}_{L^{2}}(\|\nabla u\|^2_{L^{2}}+\|\nabla\partial_{3}u\|^2_{L^{2}}).
\end{split}\tag{3.12}
\end{equation}
For $J_{2}$, we have
\begin{equation}
J_{2}\leq\frac{1}{4}\|\Delta u\|^{2}_{L^{2}}+C \|\nabla\omega\|^{2}_{L^{2}}.\tag{3.13}
\end{equation}
Substituting $(3.12),(3.13)$ in $(3.11)$, we obtain
 \begin{equation}
 \begin{split}
\frac{d}{dt}\|\nabla u\|^{2}_{L^{2}}+\|\Delta u\|^{2}_{L^{2}}
&\leq C(\|\nabla u\|^{2}_{L^{2}}+e)(\|\nabla\partial_{3}u\|^{2}_{L^{2}}+\|\nabla u\|^{2}_{L^{2}}+\|\nabla \omega\|^{2}_{L^{2}}+1).\\
\end{split}\tag{3.14}
\end{equation}
Applying the Gronwall inequality yields that
\begin{equation}
\begin{split}
\sup\limits_{0\leq t\leq T}&\|\nabla u\|^{2}_{L^{2}}\leq (e+\|\nabla u_{0}\|^{2}_{L^{2}})\exp\{\int^{T}_{0}(\|\nabla\partial_{3}u\|^{2}_{L^{2}}+\|\nabla u\|^{2}_{L^{2}}+\|\nabla \omega\|^{2}_{L^{2}}+1)d\tau\},
\end{split}\tag{3.15}
\end{equation}
which gives that
\begin{equation}
u\in L^{\infty}\big(0,T;H^{1}(\mathbb{R}^3)\big)\cap L^{2}\big(0,T;H^{2}(\mathbb{R}^3)\big).\tag{3.16}
\end{equation}
The proof of Lemma 3.1 is completed.
\medskip
\par \textbf {Proof of Theorem 1.1}\quad Taking the inner product of system $(1.1)_1$ with $-\Delta u$ in $L^2(\mathbb{R}^3)$, we obtain
\begin{equation}
\begin{split}
\frac{1}{2}\frac{d}{dt}\|\nabla u\|^{2}_{L^2}+\|\Delta u\|^{2}_{L^2}=\int_{\mathbb{R}^3}u\cdot\nabla u\cdot\Delta u dx-\int_{\mathbb{R}^3}\nabla \times\omega\cdot\Delta u dx.
\end{split}\tag{3.17}
\end{equation}
On the other hand, testing $(1.1)_2$ by $-\Delta\omega$, we get
\begin{equation}
\begin{split}
\frac{1}{2}\frac{d}{dt}\|\nabla \omega\|^{2}_{L^2}+\|\Delta \omega\|^{2}_{L^2}+2\|\nabla \omega\|^{2}_{L^2}+\|\nabla\nabla\cdot\omega\|^{2}_{L^2}
=\int_{\mathbb{R}^3}u\cdot\nabla \omega\cdot\Delta \omega dx-\int_{\mathbb{R}^3}\nabla\times u\cdot\Delta \omega dx.
\end{split}\tag{3.18}
\end{equation}
Plugging (3.17) and (3.18) together, it yields
\begin{equation}
\begin{split}
&\frac{1}{2}\frac{d}{dt}(\|\nabla u\|^{2}_{L^2}+\|\nabla \omega\|^{2}_{L^2})+\|\Delta u\|^{2}_{L^2}+\|\Delta \omega\|^{2}_{L^2}+2\|\nabla \omega\|^{2}_{L^2}+\|\nabla\nabla\cdot\omega\|^{2}_{L^2}\\
=& \int_{\mathbb{R}^3}u\cdot\nabla u\cdot\Delta u dx-\int_{\mathbb{R}^3}\nabla \times\omega\cdot\Delta u dx+\int_{\mathbb{R}^3}u\cdot\nabla \omega\cdot\Delta \omega dx-\int_{\mathbb{R}^3}\nabla\times u\cdot\Delta \omega dx\\
=&K_1+K_2+K_3+K_4.
\end{split}\tag{3.19}
\end{equation}
By using H\"{o}lder and Young inequalities, it is clear that
\begin{equation}
\begin{split}
K_{2}+K_{4}&=-\int_{\mathbb{R}^3}\nabla \times\omega\cdot\Delta u dx-\int_{\mathbb{R}^3}\nabla\times u\cdot\Delta u dx\\
&\leq C\|\nabla \times \omega\|_{L^{2}}\|\Delta u\|_{L^{2}}+C\|\nabla \times u\|_{L^{2}}\|\Delta \omega\|_{L^{2}}\\
&\leq 2C\|\nabla u\|^{2}_{L^{2}}+2C\|\nabla \omega\|^{2}_{L^{2}}+\frac{1}{2}\|\Delta u\|^{2}_{L^{2}}+\frac{1}{2}\|\Delta \omega\|^{2}_{L^{2}}.
\end{split}\tag{3.20}
\end{equation}
Combining $(3.19)$ and $(3.20)$, one obtains
\begin{equation}
\begin{split}
&\frac{d}{dt}(\|\nabla u\|^{2}_{L^2}+\|\nabla \omega\|^{2}_{L^2})+\|\Delta u\|^{2}_{L^2}+\|\Delta \omega\|^{2}_{L^2}+\|\nabla\nabla\cdot\omega\|^{2}_{L^2}\\
\leq& 4C\|\nabla u\|^{2}_{L^2}+4C\|\nabla \omega\|^{2}_{L^{2}}+2\int_{\mathbb{R}^3}u\cdot\nabla u\cdot\Delta u dx+2\int_{\mathbb{R}^3}u\cdot\nabla \omega\cdot\Delta \omega dx
\end{split}\tag{3.21}
\end{equation}
By the H\"{o}lder inequality, Gagliardo-Nirenberg inequality and Sobolev inequality,
the last two terms of $(3.21)$ can be estimated as
\begin{equation}
\begin{split}
&2\int_{\mathbb{R}^3}u\cdot\nabla u\cdot\Delta u dx+2\int_{\mathbb{R}^3}u\cdot\nabla \omega\cdot\Delta \omega dx\\
\leq& 2\|u\cdot\nabla u\|_{L^2}\|\Delta u\|_{L^2}+2\|u\cdot\nabla \omega\|_{L^2}\|\Delta \omega\|_{L^2}\\
\leq& C\|u\cdot\nabla u\|^{2}_{L^2}+C\|u\cdot\nabla \omega\|^{2}_{L^2}+\epsilon\|\Delta u\|^{2}_{L^2}+\epsilon\|\Delta \omega\|^{2}_{L^2}\\
\leq& C\| u\|^{2}_{L^\infty}(\|\nabla u\|^{2}_{L^2}+\|\nabla \omega\|^{2}_{L^2})+\epsilon\|\Delta u\|^{2}_{L^2}+\epsilon\|\Delta \omega\|^{2}_{L^2}\\
\leq& C\| \nabla u\|_{L^2}\| \Delta u\|_{L^2}(\|\nabla u\|^{2}_{L^2}+\|\nabla \omega\|^{2}_{L^2})+\epsilon\|\Delta u\|^{2}_{L^2}+\epsilon\|\Delta \omega\|^{2}_{L^2}\\
\leq& C\| \nabla u\|^{2}_{L^2}(\|\nabla u\|^{2}_{L^2}+\|\nabla \omega\|^{2}_{L^2})^2+2\epsilon\|\Delta u\|^{2}_{L^2}+\epsilon\|\Delta \omega\|^{2}_{L^2}.\\
\end{split}\tag{3.22}
\end{equation}
Plugging estimates $(3.22)$ into $(3.21)$, give us to
\begin{equation}
\begin{split}
&\frac{d}{dt}(\|\nabla u\|^{2}_{L^2}+\|\nabla \omega\|^{2}_{L^2})+\frac{1}{2}(\|\Delta u\|^{2}_{L^2}+\|\Delta \omega\|^{2}_{L^2})+\|\nabla\nabla\cdot\omega\|^{2}_{L^2}\\
\leq& 4C\|\nabla u\|^{2}_{L^2}+4C\|\nabla \omega\|^{2}_{L^2}+ C\| \nabla u\|^{2}_{L^2}(\|\nabla u\|^{2}_{L^2}+\|\nabla \omega\|^{2}_{L^2})^2.
\end{split}\tag{3.23}
\end{equation}
Setting $Y(t)=\|\nabla u\|^{2}_{L^{2}}+\|\nabla \omega\|^{2}_{L^{2}}+e$, from $(3.23)$ we see that
\begin{equation}
\frac{d}{dt}\ln Y(t)\leq 4C+ C\| \nabla u\|^{2}_{L^2}(\|\nabla u\|^{2}_{L^2}+\|\nabla \omega\|^{2}_{L^2})\ln Y(t).\tag{3.24}
\end{equation}
Due to the Gronwall inequality and Lemma 3.1, we conclude that
 \begin{equation}
\begin{split}
\ln Y(t)&\leq \bigg(\ln(e+\|\nabla u_{0}\|^{2}_{L^{2}}+\|\nabla \omega_{0}\|^{2}_{L^{2}})+4CT\bigg)\exp \int^{T}_{0}C\| \nabla u\|^{2}_{L^2}(\|\nabla u\|^{2}_{L^2}+\|\nabla \omega\|^{2}_{L^2})dt\\
&\leq \bigg(C+4CT\bigg)\exp C\int^{T}_{0}\|\nabla u\|^{2}_{L^2}+\|\nabla \omega\|^{2}_{L^2}dt.
\end{split}\tag{3.25}
\end{equation}
Combining the basic energy estimate $(3.2)$ and taking the exponential function on both side of $(3.25)$, we have
 \begin{equation}
\begin{split}
\sup\limits_{0\leq t\leq T}&(\|\nabla u\|^{2}_{L^{2}}+\|\nabla \omega\|^{2}_{L^{2}})\leq
\exp\{(C+4CT)\exp\{CT(\| u_0\|^{2}_{L^{2}}+\|\omega_0\|^{2}_{L^{2}})\}\}<\infty,
\end{split}\tag{3.26}
\end{equation}
thus
\[u,\omega\in L^{\infty}\big(0,T;H^{1}(\mathbb{R}^3)\big)\cap L^{2}\big(0,T;H^{2}(\mathbb{R}^3)\big).\]
This completes the proof of Theorem 1.1.
\section*{Acknowledgments}
\par The first author is partially supported by I.N.D.A.M-G.N.A.M.P.A. 2019 and the ``RUDN University Program 5-100''.

\end{document}